\documentclass[preprint,11pt]{elsarticle}
\usepackage[latin1]{inputenc}
\usepackage[english]{babel}
\usepackage{geometry}
\usepackage{amsmath}
\usepackage{amssymb}
\usepackage{oldgerm}
\usepackage{amscd}
\usepackage{rotating}
\usepackage[all]{xy}
\usepackage{hyperref}

\usepackage{color}

\date{\today}

\newtheorem{Prop}{Proposition}[section]
\newtheorem{Theorem}{Theorem}[section]
\newtheorem{Lemma}{Lemma}[section]
\newtheorem{Cor}{Corollary}[section]
\newcommand{\Proof}{\noindent{\bfseries Proof : }}

\newdefinition{Def}{Definition}[section]
\newdefinition{Ex}{Example}[section]
\newdefinition{Rem}{Remark}[section]

\newcommand{\Q}{\mathbb{Q}}
\newcommand{\N}{\mathbb{N}}
\newcommand{\Z}{\mathbb{Z}}
\newcommand{\IZ}{\textnormal{Int($\Z$)}}
\newcommand{\IZM}{\textnormal{Int($M_n(\Z)$)}}

\makeatletter
\def\ps@pprintTitle{%
     \let\@oddhead\@empty
     \let\@evenhead\@empty
     \def\@oddfoot{\reset@font\hfil\thepage\hfil%
     \llap{\footnotesize\itshape\today}}
     \let\@evenfoot\@oddfoot}
\makeatother

\begin{document}

\begin{frontmatter}

\title{Integral-valued polynomials over sets of algebraic integers of bounded degree}

\author{Giulio Peruginelli\fnref{fn1}}
\ead{peruginelli@math.tugraz.at}
\address{Institut f\"ur Analysis und Comput. Number Theory, Technische Universit\"at, Steyrergasse 30, A-8010 Graz, Austria.}
\fntext[fn1]{Phone number: +43 (316) 873 - 7633}

\begin{abstract}
\noindent Let $K$ be a number field of degree $n$ with ring of integers $O_K$. By means of a criterion of Gilmer for polynomially dense subsets of the ring of integers of a number field, we show that, if $h\in K[X]$ maps every element of $O_K$ of degree $n$ to an algebraic integer, then $h(X)$ is integral-valued over $O_K$, that is $h(O_K)\subset O_K$. A similar property holds if we consider the set of all algebraic integers of degree $n$ and a polynomial $f\in\Q[X]$: if $f(\alpha)$ is integral over $\Z$ for every algebraic integer $\alpha$ of degree $n$, then $f(\beta)$ is integral over $\Z$ for every algebraic integer $\beta$ of degree smaller than $n$. This second result is established by proving that the integral closure of the ring of polynomials in $\Q[X]$ which are integer-valued over the set of matrices $M_n(\Z)$ is equal to the ring of integral-valued polynomials over the set of algebraic integers of degree equal to $n$. 
\end{abstract}

\begin{keyword}
Integer-valued polynomial\sep Algebraic integers with bounded degree\sep Pr\"ufer domain \sep Polynomially dense subset \sep Integral closure \sep Pullback. MSC Classification codes: 13B25, 13F20, 11C.
\end{keyword}

\end{frontmatter}

\section{Introduction}

Let $K$ be a number field of degree $n$ over $\Q$ with ring of integers $O_K$. Given $f\in K[X]$ and $\alpha\in O_K$, the evaluation of $f(X)$ at $\alpha$ is an element of $K$. If $f(\alpha)$ is in $O_K$ we say that $f(X)$ is integral-valued on $\alpha$. If this condition holds for every $\alpha\in O_K$, we say that $f(X)$ is integral-valued over $O_K$. The set of such polynomials forms a ring, usually denoted by:
$${\rm Int}(O_K)\doteqdot\{f\in K[X]\;|\;f(O_K)\subset O_K\}.$$
Obviously, ${\rm Int}(O_K)\supset O_K[X]$ and this is a strict containment (over $\Z$, consider $X(X-1)/2$).
A classical problem regarding integral-valued polynomials is to find proper subsets $S$ of $O_K$ such that if $f(X)$ is any polynomial in $K[X]$ such that $f(s)$ is in $O_K$ for all $s$ in $S$ then $f(X)$ is integral-valued over $O_K$.  A subset $S$ of $O_K$ with this property is usually called a polynomially dense subset of $O_K$. For example, it is easy to see that cofinite subsets of $O_K$ have this property. For a general reference of polynomially dense subsets and the so-called polynomial closure see \cite{CaCh} (see also the references contained in there).  Obviously, for a polynomially dense subset $S$ we have ${\rm Int}(S,O_K)\doteqdot\{f\in K[X]\,|\,f(S)\subset O_K\}={\rm Int}(O_K)$ (in general we only have one containment). Gilmer gave a criterion which characterizes polynomially dense subsets of a 
Dedekind domain with finite residue fields (\cite{Gilmer}). His result was later elaborated by McQuillan in this way (\cite{McQ}; we state it for the ring of integers of a number field):
a subset $S$ of $O_K$ is polynomially dense in $O_K$ if and only if, for every non-zero prime ideal $P$ of $O_K$, $S$ is dense in $O_K$ with respect to the $P$-adic topology. By means of this criterion, we show here the following theorem.
\vskip0.2cm
\begin{Theorem}\label{1stthm}
Let $K$ be a number field of degree $n$ over $\Q$. Let $O_{K,n}$ be the set of algebraic integers of $K$ of degree $n$. Then $O_{K,n}$ is polynomially dense in $O_K$.
\end{Theorem}
\vskip0.2cm
The previous problem concerns the integrality of the values of a polynomial with coefficients in a number field $K$ over the set of algebraic integers of $K$. We also address here our interest to the study of the integrality of the values of a polynomial with rational coefficients over the set of algebraic integers of a proper finite extension of $\Q$,
or, more in general, over a set of algebraic integers which lie in possibly infinitely many number fields, but of bounded degree. In this direction, Loper and Werner introduced in \cite{LopWer} the following ring of integral-valued polynomials:
$${\rm Int}_{\Q}(O_K)\doteqdot\{f\in\Q[X]\;|\;f(O_K)\subset O_K\}.$$
This ring is the contraction to $\Q[X]$ of ${\rm Int}(O_K)$. It is easy to see that it is a subring of the usual ring of integer-valued polynomials $\IZ=\{f\in\Q[X]\,|\,f(\Z)\subset\Z\}$. Moreover, this is always a strict containment: take any prime integer $p$ such that there exists a prime ideal of $O_K$ above $p$ whose residue field strictly contains $\Z/p\Z$; then the polynomial $X(X-1)\ldots(X-(p-1))/p$ is in $\IZ$ but it is not in ${\rm Int}_{\Q}(O_K)$. 

In \cite{LopWer} another ring of integral-valued polynomials has been introduced. Given an integer $n$, we denote by $\mathcal{A}_n$ the set of all algebraic integers of degree less than or equal to $n$ in a fixed algebraic closure $\overline{\Q}$ of $\Q$. We define the following ring
$${\rm Int}_{\Q}(\mathcal{A}_n)\doteqdot\{f\in\Q[X]\;|\;f(\mathcal{A}_n)\subset\mathcal{A}_n\}.$$
Notice that for $n=1$ we have the usual ring $\IZ$. Given any algebraic integer $\alpha$ of degree $n$ and $f(X)$ any polynomial with rational coefficients, $f(\alpha)$ belongs to $\Q(\alpha)$ and in particular its degree is less than or equal to $n$. Therefore, a polynomial $f\in\Q[X]$ is in ${\rm Int}_{\Q}(\mathcal{A}_n)$ if and only if for every $\alpha$ in $\mathcal A_n$, $f(\alpha)$ is in the ring $\mathcal{A}_{\infty}\subset\overline{\Q}$ of all algebraic integers. We have then ${\rm Int}_{\Q}(\mathcal{A}_n)={\rm Int}_{\Q}(\mathcal{A}_n,\mathcal{A}_{\infty})=\{f\in\Q[X]\;|\;f(\mathcal{A}_n)\subset\mathcal{A}_{\infty}\}$, where the latter ring is the contraction to $\Q[X]$ of the ring ${\rm Int}(\mathcal{A}_n,\mathcal{A}_{\infty})$, which is the ring of polynomials in $\overline{\Q}[X]$ which are integer-valued over the subset $\mathcal{A}_n$ of $\mathcal{A}_\infty$ (notice that by definition the quotient field of $\mathcal{A}_\infty$ is $\overline{\Q}$). In particular, this implies that for all $n\in\N$ we have ${\rm Int}_{\Q}(\mathcal{A}_n)\subset{\rm Int}_{\Q}(\mathcal{A}_{n-1})$.  

As the authors in \cite{LopWer} show, the ring ${\rm Int}_{\Q}(\mathcal{A}_n)$ is also equal to 
\begin{equation}\label{IntQAAn}
{\rm Int}_{\Q}(\mathcal{A}_n)=\bigcap_{[K:\Q]\leq n}{\rm Int}_{\Q}(O_K)
\end{equation}
where the intersection is over the set of all number fields $K$ of degree less than or equal to $n$. One of the reason why these rings have been introduced was to show the existence of a Pr\"ufer domain lying properly between $\Z[X]$ and $\IZ$ (see question Q1 at the beginning of \cite{LopWer}). Indeed, they prove that ${\rm Int}_{\Q}(\mathcal{A}_n)$ and ${\rm Int}_{\Q}(O_K)$ for any number field $K$ are examples of such rings. 

We can ask whether a result similar to Theorem \ref{1stthm} holds for $\mathcal{A}_n$. More specifically, let $A_n$ be the subset of $\mathcal{A}_n$ of elements of degree equal to $n$. We consider the following ring
$${\rm Int}_{\Q}(A_n,\mathcal{A}_n)\doteqdot\{f\in\Q[X]\;|\;f(A_n)\subset\mathcal{A}_n\}$$
that is, polynomials with rational coefficents which map the set of algebraic integers of degree $n$ to algebraic integers (which of course have degree less than or equal to $n$; like before, notice that ${\rm Int}_{\Q}(A_n,\mathcal{A}_n)={\rm Int}_{\Q}(A_n,\mathcal{A}_{\infty})$, with the obvious notation). A priori this ring contains ${\rm Int}_{\Q}(\mathcal{A}_n)$, since $A_n$ is a subset of $\mathcal{A}_n$. By Theorem \ref{1stthm} it follows that there exist algebraic integers of degree smaller than $n$ (namely, the ring of integers of all number fields of degree $n$) on which any $f\in {\rm Int}_{\Q}(A_n,\mathcal{A}_n)$ is integral-valued (we treat this consequence of Theorem \ref{1stthm} in Remark \ref{Annotpolclo}). Using the same terminology adopted for subsets of the ring of integers of a number field, we say that $A_n$ is not polynomially closed in $\mathcal{A}_n$. The main result of this paper is that the two rings of integral-valued polynomials over $A_n $ and $\mathcal A_n$  actually coincide:
\begin{Theorem}\label{2ndthm}
For every positive integer $n$ we have: 
$${\rm Int}_{\Q}(A_n,\mathcal{A}_n)={\rm Int}_{\Q}(\mathcal{A}_n).$$
\end{Theorem}
\vskip0.4cm
In this way, in order to check whether a polynomial $f\in\Q[X]$ is integral-valued over the set of algebraic integers of degree bounded by $n$, it is sufficient to check if $f(X)$ is integral-valued over the set of algebraic integers of degree exactly equal to $n$. As above, we say that $A_n$ is polynomially dense in $\mathcal{A}_n$. In particular, Theorem \ref{1stthm} and Theorem \ref{2ndthm} also show that it is not necessary to take a set $S$ of algebraic integers properly containing $\Z$ to exhibit a Pr\"ufer domain properly contained between $\Z[X]$ and $\IZ$ (specifically, ${\rm Int}_{\Q}(S,\mathcal{A}_{\infty})$). 

The ring ${\rm Int}_{\Q}(\mathcal{A}_n)$ has been introduced in \cite{LopWer} also for another reason. We denote by $M_n(\Z)$ the noncommutative $\Z$-algebra of $n\times n$ matrices with entries in $\Z$. Given a matrix $M\in M_n(\Z)$ and a polynomial $f\in\Q[X]$, $f(M)$ is a matrix with rational entries. If $f(M)$ has integer entries, we say that $f(X)$ is integer-valued on $M$. We consider the ring of polynomials which are integer-valued over all the matrices of $M_n(\Z)$ (introduced in \cite{Frisch0}):
$$\IZM\doteqdot\{f\in\Q[X]\;|\;f(M)\in M_n(\Z),\,\forall M\in M_n(\Z)\}.$$
Like ${\rm Int}_{\Q}(\mathcal{A}_n)$, $\IZM$ is a subring of the ring of integer-valued polynomial $\IZ$ ($\Z$ embeds into $M_n(\Z)$ as the subalgebra of scalar matrices). In \cite{LopWer} the following theorem is proved:
\vskip0.3cm
\begin{Theorem}\label{LW}
The ring $\IZM$ is not integrally closed and its integral closure is ${\rm Int}_{\Q}(\mathcal{A}_n)$, which is a Pr\"ufer domain.
\end{Theorem}
What is behind the containment between these two rings of integer-valued polynomials is based on the following fact: if $f\in\Q[X]$ is integer-valued over a matrix $M\in M_n(\Z)$, then $f(X)$ is integral-valued over the set of roots of the characteristic polynomial of $M$, which are elements of $\mathcal{A}_n$ (that is, the eigenvalues of $M$ as a matrix over $\overline{\Q}$). This fact can be easily proved directly, but we will give an alternative point of view based on a representation of the ring $\IZM$ as an intersection of pullbacks of $\Q[X]$, which turn out to be equal to a particular class of rings of integer-valued polynomials over algebraic integers. Our Theorem \ref{2ndthm} is proved by means of a result of Frisch (\cite[Proposition 6.2]{Frisch1}, \cite{Frischcor}) which says that the set of matrices with irreducible characteristic polynomial is polynomially dense in $M_n(\Z)$.

We give a summary of the paper. In the second section we use Gilmer's criterion to prove a generalization of  Theorem \ref{1stthm} for orders of a number field $K$. By means of a result of Gy{\H{o}}ry, we also show that if we remove from $O_{K,n}$ those elements which generates $O_K$ as a $\Z$-algebra (if any) we still get a polynomially dense subset. 

In Section \ref{RingIntPullback} we establish the connection between the subring of $\Q[X]$ of integer-valued polynomials over a set of $n\times n$ integral matrices with irreducible characteristic polynomial in a prescribed set $P$ (these rings have been studied in \cite{Per2}) and rings of integral-valued polynomials over the roots of the polynomials in $P$ (which are algebraic integers of degree $n$), showing that the latter is the integral closure of the former. This result is achieved by means of Theorem \ref{IntCloSARA}. We need a result from \cite{Per2} which says that the ring of integer-valued polynomials over matrices with prescribed characteristic polynomial $p\in\Z[X]$ is equal to the pullback ring $\Z[X]+p(X)\cdot \Q[X]$. Theorem \ref{IntCloSARA} proves that the integral closure of the intersection of a collection of pullbacks of this kind is the intersection of the integral closure of the pullbacks of the family. The crucial assumption we need is that the degree of the polynomials of the set $P$ is bounded by the above integer $n$. From this result, Theorem \ref{2ndthm} follows. In Remark \ref{Annotpolclo} we stress that Theorem \ref{1stthm} does not imply directly Theorem \ref{2ndthm}.

We also prove another related result, analogous to Theorem \ref{LW}. As ${\rm Int}_{\Q}(\mathcal{A}_n)$ is the integral closure of the ring ${\rm Int}(M_n(\Z))$, in the same way, for every number field $K$, ${\rm Int}_{\Q}(O_K)$ (which by (\ref{IntQAAn}) is an overring of ${\rm Int}_{\Q}(\mathcal{A}_n)$) is equal to the integral closure of a certain overring of ${\rm Int}(M_n(\Z))$, namely the ring of integer-valued polynomials over the set $M_n^K(\Z)$ of matrices with characteristic polynomial equal to a minimal polynomial of some algebraic integer of maximal degree of $K$.

\vskip1cm

\section{Polynomially dense subsets of the ring of integers of a number field}\label{poldensesubsetOK}
\vskip0.4cm
We recall from the introduction that a subset $S$ of a domain $D$ with quotient field $K$ is polynomially dense in $D$ if the ring  
${\rm Int}(S,D)\doteqdot\{f\in K[X]\,|\,f(S)\subset D\}$ is equal to ${\rm Int}(D)=\{f\in K[X]\,|\,f(D)\subset D\}$. 

Let $K$ be a number field with ring of integers $O_K$. The statement of the criterion of Gilmer (\cite[Theorem 8]{Gilmer}, originally stated for Dedekind domains with finite residue fields) we mentioned in the introduction is the following: a subset $S$ of $O_K$ is polynomially dense in $O_K$ if and only if for every non-zero prime ideal $P$ of $O_K$ and any positive integer $k$, the set $S$ contains a complete set of representatives of the residue classes modulo $P^k$. According to Gilmer's terminology, such a subset is called prime power complete. As McQuillan showed in \cite{McQ}, this property corresponds to $S$ being dense in $O_K$ with respect to the $P$-adic topology, for each non-zero prime ideal $P$ of $O_K$. By \cite[Chapter IV]{CaCh} the same result holds for orders of $K$. We recall that an order of $K$ is a subring of $K$ which has maximal rank as a  $\Z$-module; in particular, an order $O$ of $K$ is contained in $O_K$, which is called maximal order. If $n$ is the degree of $K$ over $\Q$ and $O\subseteq O_K$ is an order, we denote by $O_{n}$ the set of elements of $O$ of degree $n$, thus $O_n\doteqdot\{\alpha\in O\,|\,\Q(\alpha)=K\}$.  We denote ${\rm Int}(O_{n},O)=\{f\in K[X]\,|\,f(O_n)\subset O\}$ and ${\rm Int}(O)=\{f\in K[X]\,|\,f(O)\subset O\}$ (the quotient field of $O$ is the number field $K$). Via Gilmer's criterion it is easy to see that, if a subset $S$ of $O_K$ is polynomially dense in $O_K$ then it has non-zero intersection with $O_n$. In fact, if $S$ is contained in a proper subfield $K'$ of $K$, just consider a prime ideal $P$ of $O_K$ whose residue field is strictly bigger than the residue field of the contraction of $P$ to $O_{K'}$. This implies that there are residue classes modulo $P$ which are not covered by the set $S$.

The next theorem is a generalization of Theorem \ref{1stthm} of the introduction.

\begin{Theorem}\label{poldense}
Let $K$ be a number field of degree $n$ and $O$ an order of $K$. Then
$${\rm Int}(O_{n},O)={\rm Int}(O)$$
that is, $O_{n}$ is polynomially dense in $O$.
\end{Theorem}

\Proof  Given a non-zero ideal $I$ of $O$, there exists $\alpha\in I$ of degree $n$. In fact, suppose $I\cap\Z=d\Z$, for some non-zero integer $d$. Since $dO\subset I$, it is sufficient to prove the claim for the principal ideals of $O$ generated by an integer $d$. Pick an element $\alpha$ in $O$ of degree $n$. Then there are exactly $n$ conjugates of $d\alpha$ over $\Q$ and they lie in $dO$. 

Let $P^k$ be the power of a non-zero prime ideal of $O_K$. By the discussion above $P^k$ has non-trivial intersection with $O_{n}$. We have to show that every residue class $\alpha+P^k$ has a representative which lies in $O_{n}$. If $\alpha\in O\setminus O_{n}$, pick an element $\beta\in P^k$ of maximal degree. Then, using an argument similar to the proof of the primitive element theorem, there is an integer $k$ such that the algebraic integer $\gamma=\alpha+k\beta$ (which is in $\alpha+P^k$) is a generator of $\Q(\alpha,\beta)=\Q(\beta)=K$, thus $\gamma$ has maximal degree $n$. $\Box$ 

\vskip0.5cm

Since $O_{K,n}\subset O_K\setminus\Z\subset O_K$ we also have that $O_K\setminus\Z$ is polynomially dense in $O_K$.

We may wonder whether there exist polynomially dense subsets properly contained in the subset $O_n$ of an order $O$ of a number field $K$ of degree $n$. The next proposition gives a positive answer to this question. Given such an order $O\subseteq O_K$, we consider the set $$A_O\doteqdot\{\alpha\in O\,|\,\Z[\alpha]=O\}.$$
The set $A_O$ is contained in $O_n$ and it may be empty. By a result of Gy{\H{o}}ry, 
 $A_O$ is a finite union of equivalence classes with respect to the equivalence relation given by $\alpha\sim\beta \Leftrightarrow \beta=\pm\alpha+m$, for some $m\in\Z$ (see \cite{Gyory}). This means that
$$A_O=\bigcup_{i=1,\ldots,k}(\pm\alpha_i+\Z)$$
where $O=\Z[\alpha_i]$, $i=1,\ldots,k$ and $\alpha_i\pm\alpha_j\notin\Z$, $\forall i\not=j$. With these notations we have the following Proposition.
\vskip0.6cm
\begin{Prop}\label{OnApoldenseO}
$O_{n}\setminus A_O$ is polynomially dense in $O$.
\end{Prop}

\Proof If $A_O$ is empty then by Theorem \ref{poldense} we are done. Suppose now that $A_O$ is not empty. Let $I=P^k$ be a power of a prime ideal $P$ (our arguments hold for any ideal $I$ of $O$, indeed). Suppose that a residue class $\alpha+I$ is contained in $A_O$. This means that the class $\alpha+I$ itself (and consequently the ideal $I$) can be partitioned into a finite union of sets, each contained in $\pm\beta_i+\Z$,  $\beta_i=\alpha_i+\alpha$. That is, we have $I=\bigcup_{1,\ldots,k}(\pm\beta_i+\Z)\cap I$. 
Now for each $i=1,\ldots,k$, choose (if it exists) $\gamma_i\in I$ such that $\gamma_i-\beta_i\in \Z$ (there exists at least one such value of $i$, or else $I$ would be empty). Then we have
$$I=\bigcup_{i=1,\ldots,t}\left(\gamma_i+(\Z\cap I)\right)$$
($t\leq k$; the containment $(\supseteq)$ is obvious; 
conversely, if $\gamma\in I$, for some $\beta_i$ we have $\beta_i\sim\gamma$,
 so that $\gamma\sim\gamma_i$ and so $\gamma-\gamma_i\in\Z\cap I$). Hence, the additive group of $I$ is a finite union of residue classes modulo
$J=I\cap\Z$. This is not possible: $J$ is a free-$\Z$ module of rank $1$ and $I$ is a free-$\Z$ module of
rank $n>1$.
$\Box$

\vskip1cm

\section{Rings of integer-valued polynomials as intersection of pullbacks}\label{RingIntPullback}
\vskip0.4cm
Let $\alpha$ be a fixed algebraic integer over $\Z$. Given a polynomial $f\in\Q[X]$, the evaluation of $f(X)$ at $\alpha$ is an element of the number field $K=\Q(\alpha)$, which is by its very definition the set of all the $f(\alpha)$'s, with $f\in\Q[X]$. This set is clearly a ring (a field, indeed), since it is the image under the evaluation homomorphism at $\alpha$ of the polynomial ring $\Q[X]$. It is well-known that the set of those $f(\alpha)$'s which are integral over $\Z$ is a subring of $K$, called the ring of integers of $K$. We denote this subring by $O_{\Q(\alpha)}$ (we stress that this subring does not depend on $\alpha$ but only on the number field $K$). Notice that, if $f\in\Z[X]$, then $f(\alpha)$ is clearly in $\Z[\alpha]\subseteq O_{\Q(\alpha)}$, but there exist other polynomials $f(X)$ in $\Q[X]\setminus\Z[X]$ such that $f(\alpha)$ is in $\Z[\alpha]$. In order to study this phenomenon, we introduce the following rings.
\begin{Def}
\begin{align*}
R_\alpha\doteqdot{\rm Int}_{\Q}(\{\alpha\},\Z[\alpha])=&\{f\in\Q[X]\,|\,f(\alpha)\in\Z[\alpha]\}  \\
S_\alpha\doteqdot{\rm Int}_{\Q}(\{\alpha\},O_{\Q(\alpha)})=&\{f\in\Q[X]\,|\,f(\alpha)\in O_{\Q(\alpha)}\}.
\end{align*}
\end{Def}
Notice that $\Z[X]\subset R_\alpha\subseteq S_\alpha\subset\Q[X]$, so that $R_\alpha$ and $S_\alpha$ are $\Z[X]$-algebras. 

If $\alpha$ generates the ring $O_{\Q(\alpha)}$ as a $\Z$-algebra, that is $\Z[\alpha]=O_{\Q(\alpha)}$, we say that the ring of integers $O_{\Q(\alpha)}$ is monogenic. In particular, if this condition holds, $R_\alpha=S_\alpha$. It is easy to see that in general the containment $R_\alpha\subset S_\alpha$ is proper.
Take $\alpha=2\sqrt{2}$. Then $\Q(\alpha)=\Q(\sqrt{2})$ so that $O_{\Q(\alpha)}=\Z[\sqrt{2}]$. We consider $f(X)=X/2$. Then $f(\alpha)=\sqrt{2}\in O_{\Q(\alpha)}\setminus \Z[\alpha]$, so that $f\in S_\alpha\setminus R_\alpha$. In general $\Z[\alpha]$ is only contained in $O_{\Q(\alpha)}$ and the two rings have the same quotient field $\Q(\alpha)$. The integral closure of $\Z[\alpha]$ in $\Q(\alpha)$ is obviously $O_{\Q(\alpha)}$. By the next lemma, the previous implication can be reversed, namely if $R_\alpha=S_\alpha$ then $\Z[\alpha]=O_{\Q(\alpha)}$. 

\begin{Lemma}\label{RalphaSalpha}
Let $\alpha$ be an algebraic integer and $K=\Q(\alpha)$. We set
$$
R_{\alpha}(\alpha)\doteqdot\{f(\alpha)\,|\,f\in R_\alpha\},\hskip0.5cm
S_{\alpha}(\alpha)\doteqdot\{f(\alpha)\,|\,f\in S_\alpha\}.
$$
Then $R_{\alpha}(\alpha)=\Z[\alpha]$ and $S_{\alpha}(\alpha)=O_K$.
\end{Lemma}
\Proof  By definition $R_{\alpha}(\alpha)\subseteq\Z[\alpha]$. Since $\Z[X]\subset R_{\alpha}$ we also have the other containment.
For the same reason we have $S_{\alpha}(\alpha)\subseteq O_K$. Let $c=c_\alpha=[O_K:\Z[\alpha]]$ and take $\beta\in O_K$.  We have that $c\beta\in\Z[\alpha]$, so that $c\beta=g(\alpha)$ for some $g\in\Z[X]$. Then $f(X)\doteqdot g(X)/c\in\Q[X]$ has the property that $f(\alpha)=\beta\in O_K$, that is $f\in S_{\alpha}$ and its evaluation on $\alpha$ is $\beta$ as wanted. $\Box$
\vskip0.3cm

By \cite[Proposition IV.4.3]{CaCh}, the integral closure of ${\rm Int}(\{\alpha\},\Z[\alpha])$ in its quotient field $\Q(\alpha)(X)$ is ${\rm Int}(\{\alpha\},O_K)$. We notice that ${\rm Int}(\{\alpha\},O_K)=O_K+(X-\alpha)K[X]$, where $K=\Q(\alpha)$, so ${\rm Int}(\{\alpha\},O_K)$ is a pullback of $K[X]$. The next proposition shows that analogous properties hold for the contraction of these rings to $\Q[X]$, which are the rings  $R_\alpha$ and $S_\alpha$, respectively. For a general treatment of pullbacks see \cite{GabHou}.
\vskip0.5cm
\begin{Prop}\label{RSalphapullbacks}
Let $\alpha$ be an algebraic integer of degree $n$ over $\Z$ and $p_\alpha\in\Z[X]$ its minimal polynomial. Then $R_\alpha$ and $S_\alpha$ are pullbacks of $\Q[X]$. In particular, $R_\alpha=\Z[X]+M_\alpha$, where $M_\alpha$ is the maximal ideal of $\Q[X]$ generated by $p_\alpha(X)$. Moreover, the integral closure of $R_\alpha$ in its quotient field $\Q(X)$ is $S_\alpha$, which is a Pr\"ufer domain, and $R_\alpha$ is integrally closed if and only if $\Z[\alpha]=O_{\Q(\alpha)}$.
\end{Prop}
\vskip0.2cm
\Proof It is easy to see that $M_\alpha=p_{\alpha}(X)\cdot\Q[X]$ is a common ideal of $R_\alpha$, $S_\alpha$ and $\Q[X]$. We have the following diagram:
\begin{equation}\label{*}
\begin{array}{ccc}
R_\alpha&\to&R_\alpha/M_\alpha\cong\Z[\alpha]\\
\downarrow&&\downarrow\\
S_\alpha&\to&S_\alpha/M_\alpha\cong O_{\Q(\alpha)}\\
\downarrow&&\downarrow\\
\Q[X]&\to&\Q[X]/M_\alpha\cong \Q(\alpha)
\end{array}
\end{equation}
where the vertical arrows are the natural inclusions and the horizontal arrows are the natural projection, which can be viewed as the evaluation of a polynomial $f(X)$ at $\alpha$ (in fact, the residue class $f(X)+M_\alpha$ is equal to $f(\alpha)+M_\alpha$). The kernel of the evaluation homomorphism at $\alpha$ at each row is the ideal $M_\alpha$. Notice that $M_\alpha\cap\Z=\{0\}$, so that $\Z$ injects into the above residue rings. Because of that and by Lemma \ref{RalphaSalpha}, $R_\alpha/M_\alpha\cong R_\alpha(\alpha)=\Z[\alpha]$ and $S_\alpha/M_\alpha\cong S_\alpha(\alpha)=O_K$. Obviously, $\Q[X]/M_\alpha\cong\Q(\alpha)$. Thus, by definition, $R_\alpha$ and $S_\alpha$ are pullbacks of $\Q[X]$. 

We have to show that $R_\alpha=\Z[X]+M_\alpha$. The containment $(\supseteq)$ is straightforward. Conversely, let $f\in R_\alpha$, $f(X)=g(X)/d$, for some $g\in\Z[X]$ and $d\in\Z\setminus\{0\}$. Then $g(X)=q(X)p_\alpha(X)+r(X)$, for some $q,r\in\Z[X]$, $\deg(r)<n$. Then $f(\alpha)=r(\alpha)/d\in\Z[\alpha]$, that is $r(\alpha)\in d\Z[\alpha]$. Now $1,\alpha,\ldots,\alpha^{n-1}$ is a free $\Z$-basis of the $\Z$-module $\Z[\alpha]$. This means that, if  $r(\alpha)=\sum_{i=0,\ldots,n-1}a_i\alpha^i$, for some $a_0,\ldots,a_{n-1}\in \Z$, then $d$ divides $a_i$ for all $i$, so that $r\in d\Z[X]$. This shows that $f\in \Z[X]+M_\alpha$.

Finally, since $O_{\Q(\alpha)}$ is the integral closure of $\Z[\alpha]$, we apply \cite[Theorem 1.2]{GabHou} to our pullback diagram to conclude that $S_\alpha$ is the integral closure of $R_\alpha$ and $R_\alpha$ is integrally closed if and only if $\Z[\alpha]=O_{\Q(\alpha)}$. Finally, in a pullback diagram like in (\ref{*}) (see for instance \cite[Corollary 4.2]{GabHou}), $S_\alpha$ is a Pr\"ufer domain because $\Q[X]$ and $O_{\Q(\alpha)}$ are Pr\"ufer domains (indeed, they are Dedekind domains). Notice that, for the same result, if $\Z[\alpha]\subsetneq O_K$, $R_\alpha$ cannot be Pr\"ufer, since in this case $\Z[\alpha]$ is not integrally closed.  $\Box$
\vskip0.5cm

We now consider an arbitrary set $\mathcal{A}$ of algebraic integers of bounded degree and the corresponding intersections of the rings $R_\alpha$ and $S_\alpha$, for $\alpha\in\mathcal{A}$, respectively. We recall from the introduction that for each positive integer $n$, $\mathcal{A}_n$ denotes the set of algebraic integers of degree bounded by $n$ and $A_n$ the subset of $\mathcal{A}_n$ of those algebraic integers of degree equal to $n$.

\begin{Def} Given a subset $\mathcal{A}$ of $\mathcal{A}_n$, we set
$$\mathcal{R}_\mathcal{A}\doteqdot\bigcap_{\alpha\in\mathcal{A}}R_\alpha\subseteq\mathcal{S}_\mathcal{A}\doteqdot\bigcap_{\alpha\in\mathcal{A}}S_\alpha.$$
\end{Def}

\begin{Rem}\label{PropRASA}
Notice that, for all $\mathcal{A}\subseteq\mathcal{A}_n$ we have $\Z[X]\subset \mathcal{R}_\mathcal{A}\subseteq\mathcal{S}_\mathcal{A}\subset\Q[X]$, so that $\mathcal{R}_\mathcal{A}$ and $\mathcal{S}_\mathcal{A}$ are $\Z[X]$-algebras. For every $\mathcal{A}\subseteq\mathcal{A}_n$, we have $\mathcal{S}_{\mathcal{A}_n}\subseteq\mathcal{S}_\mathcal{A}$. This containment implies that $\mathcal{S}_\mathcal{A}$ is Pr\"ufer, since it is an overring of the Pr\"ufer ring ${\rm Int}_{\Q}(\mathcal{A}_n)$ (Theorem \ref{LW}). 

Let $A_\alpha\doteqdot\{\alpha_1=\alpha,\ldots,\alpha_n\}$ be the set of conjugates of $\alpha$ over $\Q$. If we consider the action of the Galois group of $\Q(\alpha_1,\ldots,\alpha_n)$ over $\Q$, it is easy to see that $R_{\alpha_i}=R_{\alpha_j}, S_{\alpha_i}=S_{\alpha_j}$ for all $i,j\in\{1,\ldots,n\}$. In particular, $S_\alpha=\mathcal{S}_{A_\alpha}$ and $R_\alpha=\mathcal{R}_{A_\alpha}$. Given an algebraic integer $\alpha$ of degree $m\leq n$, it is easy to see that the following equality holds: 
$$S_\alpha={\rm Int}_{\Q}(\{\alpha\},\mathcal{A}_n)\doteqdot\{f\in\Q[X]\,|\,f(\alpha)\in\mathcal{A}_n\}.$$
It follows that, for any $\mathcal{A}\subseteq \mathcal{A}_n$, we have $\mathcal{S}_\mathcal{A}={\rm Int}_{\Q}(\mathcal{A},\mathcal{A}_n)=\{f\in\Q[X]\,|\,f(\mathcal A)\subset\mathcal{A}_n\}$. In the same way, as we remarked in the introduction, $\mathcal{S}_\mathcal{A}={\rm Int}_{\Q}(\mathcal{A},\mathcal{A}_\infty)$. From this fact it follows that the ring $\mathcal{S}_\mathcal{A}$ is integrally closed by \cite[Proposition IV.4.1]{CaCh}, since $\mathcal{S}_\mathcal{A}$ is the contraction to $\Q(X)$ of an integrally closed subring of $\overline{\Q}(X)$. 

In particular, we have
\begin{equation}\label{SAnAn}
\mathcal{S}_{\mathcal{A}_n}={\rm Int}_{\Q}(\mathcal{A}_n)\subseteq {\rm Int}_{\Q}(A_n,\mathcal{A}_n)=\mathcal{S}_{A_n}.
\end{equation}

This is another representation of the rings ${\rm Int}_{\Q}(\mathcal{A}_n)$ and ${\rm Int}_{\Q}(A_n,\mathcal{A}_n)$ as an intersection of the rings $S_\alpha$. Notice that $\mathcal{S}_\Z=\IZ$. In the same way, given a number field $K$, $[K:\Q]\leq n$,  we have 

$${\rm Int}_{\Q}(O_K)={\rm Int}_{\Q}(O_K,\mathcal{A}_n)=\mathcal S_{O_K}.$$

So for all number fields $K$, the ring ${\rm Int}_{\Q}(O_K)$ can be represented as an intersection of the rings $S_\alpha$, $\alpha\in O_K$. Actually, by Theorem \ref{poldense} we can restrict the intersection to the algebraic integers of $K$ of degree $n=[K:\Q]$: $\mathcal S_{O_K}=\mathcal S_{O_{K,n}}$. In most cases, we can even consider the intersection on strictly smaller subsets of algebraic integers which are polynomially dense in $O_K$ (see Proposition \ref{OnApoldenseO}).
\end{Rem} 

\begin{Rem}\label{Annotpolclo} Since $A_n=\bigcup_{[K:\Q]=n}O_{K,n}$ (the union is over the family of number fields of degree $n$ over $\Q$) we have $\mathcal{S}_{A_n}=\bigcap_{[K:\Q]=n}\mathcal{S}_{O_{K,n}}$ and $\mathcal{R}_{A_n}=\bigcap_{[K:\Q]=n}\mathcal{R}_{O_{K,n}}$. By the above observation that $\mathcal{S}_{O_{K,n}}=\mathcal{S}_{O_{K}}$, we have
\begin{equation}\label{IntAn}
{\rm Int}_{\Q}(A_n,\mathcal{A}_n)=\bigcap_{[K:\Q]=n}{\rm Int}_{\Q}(O_K).
\end{equation}
This representation for ${\rm Int}_{\Q}(A_n,\mathcal{A}_n)$ is similar to the representation for ${\rm Int}_{\Q}(\mathcal{A}_n)$ in (\ref{IntQAAn}) (notice that this intersection is over all the number fields $K$ of degree equal to $n$, while in (\ref{IntQAAn}) the intersection is taken over all the number fields of degree less than or equal to $n$). In particular, the equality in (\ref{IntAn}) shows that $A_n$ is not polynomially closed in $\mathcal{A}_n$. More precisely, given a polynomial $f\in\Q[X]$ which is integral-valued over the set of all the algebraic integers of degree $n$, it follows that $f(X)$ is integral-valued over the ring of integers of every number field of degree $n$. We notice that from (\ref{IntAn}) we have Theorem \ref{2ndthm} for $n=2$. However, this not prove Theorem \ref{2ndthm} in general, because the algebraic integers of a number field of degree $n$ have degree which divides $n$ (for example, if $n=3$ then $O_K=O_{K,3}\cup\Z$, no algebraic integers of degree $2$ can be in $O_K$). 
\end{Rem}
\vskip0.3cm
We give now a generalization of the last statement of Proposition \ref{RSalphapullbacks}. The next theorem shows that, given any subset  $\mathcal{A}$ of algebraic integers of degree bounded by $n$, the integral closure of $\mathcal{R}_\mathcal{A}$ is $\mathcal{S}_\mathcal{A}$, so, by Proposition \ref{RSalphapullbacks}, we can say that the integral closure of the intersection of the rings $R_{\alpha}$, $\alpha$ in $\mathcal{A}$, is equal to the intersection of their integral closures $S_{\alpha}$.

\begin{Theorem}\label{IntCloSARA}
For any $\mathcal{A}\subseteq\mathcal{A}_n$, $\mathcal{S}_{\mathcal{A}}$ is the integral closure of $\mathcal{R}_{\mathcal{A}}$.
\end{Theorem}
\vskip0.2cm
The proof of Theorem \ref{IntCloSARA} follows by the argument given in \cite{LopWer} to show that the integral closure of ${\rm Int}(M_n(\mathbb{Z}))$ is ${\rm Int}_{\Q}(\mathcal{A}_n)$. 
\vskip0.6cm
\begin{Lemma}\label{almintSARA}
For all $f\in \mathcal{S}_{\mathcal{A}}$, there exists $c\in\Z\setminus\{0\}$ such that $c\cdot\mathcal{R}_\mathcal{A}[f]\subset\mathcal{R}_\mathcal{A}$.
\end{Lemma}
\vskip0.2cm
The lemma says that every element of $\mathcal{S}_{\mathcal{A}}$ is almost integral over $\mathcal{R}_{\mathcal{A}}$, that is, $\mathcal{S}_{\mathcal{A}}$ is contained in the complete integral closure of $\mathcal{R}_{\mathcal{A}}$ (remember that both have the same quotient field $\Q(X)$). In particular, we also have $c\cdot \Z[f]\subset\mathcal{R}_{\mathcal{A}}$.
\vskip0.2cm
\Proof It is sufficient to show that there exists a non-zero $c\in\Z$ such that for every $i\in\N$, $c\cdot f(X)^i\in \mathcal{R}_\mathcal{A}$. Let $i\in\N$ be fixed and let $\alpha\in\mathcal{A}$. We know that $f(\alpha)\in O_{\Q(\alpha)}$, so there exists a monic polynomial $m_\alpha\in\Z[X]$ of degree $\leq n$ such that $m_\alpha(f(\alpha))=0$. For all $\alpha\in\mathcal{A}$ we have
$$X^i=q_{\alpha,i}(X)m_\alpha(X)+r_{\alpha,i}(X)$$
for some $q_{\alpha,i},r_{\alpha,i}\in\Z[X]$, $r_{\alpha,i}(X)=0$ or $\deg(r_{\alpha,i})<n$.
Therefore
$$f(\alpha)^i=r_{\alpha,i}(f(\alpha)).$$
Since there is a uniform bound on the degree of the polynomials $r_{\alpha,i}(X)$, for $\alpha\in\mathcal{A}$ and $i\in\N$, $c\cdot f(\alpha)^i\in\Z[\alpha]$ for some $c\in\Z$ (actually, we can take $c=d^{n-1}$, where $d$ is a common denominator of the coefficients of $f(X)$). Notice that $c$ does not depend on $i$. It follows that $c\cdot f(X)^i\in\mathcal{R}_\mathcal{A}$ for all $i\in\N$. $\Box$
\vskip0.3cm
The following proposition will prove Theorem \ref{IntCloSARA}, since $\mathcal{S}_\mathcal{A}$ is integrally closed by Remark \ref{PropRASA}.
\vskip0.3cm
\begin{Prop}
For every $\mathcal{A}\subseteq\mathcal{A}_n$, $\mathcal{R}_\mathcal{A}\subseteq\mathcal{S}_\mathcal{A}$ is an integral ring extension.
\end{Prop}
\Proof Let $f\in \mathcal{S}_\mathcal{A}$ and let $c\in\Z$ as in Lemma \ref{almintSARA}. We show that there exists a monic $\varphi\in\Z[X]$ such that $\varphi(f(X))\in\mathcal{R}_\mathcal{A}$. 

By hypothesis, for each $\alpha\in\mathcal{A}$ there exists $m_\alpha\in\Z[X]$ monic of degree $n$ such that $m_\alpha(f(\alpha))=0$. Let $S$ be a set of monic representatives in $\Z[X]$ of all the degree $n$ monic polynomials in the quotient ring $(\Z/c^2\Z)[X]$ and let $\varphi(X)=\prod_{\varphi_i\in S}\varphi_i(X)\in\Z[X]$ (notice that $\varphi$ is monic). For each $\alpha\in\mathcal{A}$ there exists $i=i(\alpha)$ such that $m_\alpha(X)\equiv\varphi_i(X)\pmod{c^2}$. Let $\varphi=\varphi_i\widetilde{\varphi}$ where $\widetilde{\varphi}$ is the product of the remaining $\varphi_j$ in $S$. In this way we have
$$\varphi(f(\alpha))=c^2 r(f(\alpha))\widetilde{\varphi}(f(\alpha))$$
for some $r\in\Z[X]$. By the previous lemma $\varphi(f(\alpha))$ is in $\Z[\alpha]$, and this holds for any $\alpha\in\mathcal{A}$ ($\varphi$ is independent of $\alpha$). Hence, $\varphi(f(X))\in \mathcal{R}_\mathcal{A}$. $\Box$
\vskip0.5cm
\begin{Rem}\label{RAintclo} We conclude this section with a remark. The ring $\mathcal{R}_\mathcal{A}$ is integrally closed if and only if it is a Pr\"ufer domain, because its integral closure $\mathcal{S}_\mathcal{A}$ is always a Pr\"ufer domain. If this holds, then all the overrings $R_\alpha$ of $\mathcal{R}_\mathcal{A}$, for $\alpha\in\mathcal{A}$, are integrally closed. By Proposition \ref{RSalphapullbacks} this holds if and only if $\Z[\alpha]=O_{\Q(\alpha)}$ for each such $\alpha$'s. If the last condition holds, clearly $\mathcal{R}_\mathcal{A}=\mathcal{S}_\mathcal{A}$.  We have thus shown that
\begin{center}
$\mathcal{R}_\mathcal{A}$ is integrally closed if and only if $\mathcal{A}\subseteq \widehat{\mathcal A}_n,$
\end{center}
where $\widehat{\mathcal A}_n\doteqdot\{\alpha\in\mathcal{A}_n\,|\,\Z[\alpha]=O_{\Q(\alpha)}\}$. If $\mathcal{A}\subseteq O_{K,n}$, where $K$ is a given number field, then $\mathcal{R}_\mathcal{A}=\mathcal{S}_\mathcal{A}$ if and only if $\mathcal{A}\subseteq A_{O_K}=\{\alpha\in O_K\,|\,\Z[\alpha]=O_K\}$.
\end{Rem}
\vskip0.3cm
The next subsection gives a connection between ring of integral-valued polynomials over algebraic integers and rings of integer-valued polynomials over integral matrices.

\subsection{\textbf{Integer-valued polynomials over matrices}}

We recall some results from \cite{Per2}. Given a subset $\mathcal{M}$ of $M_n(\Z)$, we consider the ring of polynomials which are integer-valued over $\mathcal{M}$:
$${\rm Int}(\mathcal{M},M_n(\Z))\doteqdot\{f\in \Q[X]\;|\;f(M)\in M_n(\Z),\forall M\in\mathcal{M}\}.$$
Let $\mathcal{P}_n$ be the set of monic polynomials in $\Z[X]$ of degree $n$. For a subset $P$ of $\mathcal{P}_n$ we denote by $M_n^P(\Z)$ the subset of $M_n(\Z)$ of those matrices whose characteristic polynomial is in $P$. If $P=\{p(X)\}$, we set $M_n^p(\Z)=M_n^P(\Z)$; this is the set of matrices whose characteristic polynomial is equal to $p(X)$. For a generic $P\subseteq\mathcal{P}_n$, since $M_n^P(\Z)=\bigcup_{p\in P}M_n^p(\Z)$ we have 
\begin{equation}\label{critMnP}
{\rm Int}(M_n^P(\Z),M_n(\Z))=\bigcap_{p\in P}{\rm Int}(M_n^p(\Z),M_n(\Z)).
\end{equation}
We recall that, by \cite[Lemma 2.2 and Remark 2.1]{Per2}, for every $p\in\mathcal{P}_n$ we have
\begin{equation}\label{IntMnppull}
{\rm Int}(M_n^p(\Z),M_n(\Z))=\Z[X]+p(X)\cdot \Q[X].
\end{equation}
In this way, by (\ref{critMnP}) above, every ring of integer-valued polynomial over matrices with prescribed characteristic polynomial can be represented as an intersection of pullbacks of $\Q[X]$ of the above kind. By Proposition \ref{RSalphapullbacks} and (\ref{IntMnppull}) we have the following corollary, which establishes the connection between rings of integer-valued polynomials over matrices and rings of integral-valued polynomials over algebraic integers:
\vskip0.4cm
\begin{Cor}\label{IntalphaZalpha=IntMnpalpha}
Let $\alpha$ be an algebraic integer of degree $n$ over $\Z$ and $p_\alpha\in\Z[X]$ its minimal polynomial. Then
$${\rm Int}_{\Q}(\{\alpha\},\Z[\alpha])={\rm Int}(M_n^{p_{\alpha}}(\Z),M_n(\Z))=\Z[X]+p_{\alpha}(X)\cdot \Q[X].$$
\end{Cor}

\begin{Rem} Immediately from (\ref{IntMnppull}) and by the Cayley-Hamilton theorem, for every $M\in M_n^p(\Z)$ the evaluation of any $f\in {\rm Int}(M_n^p(\Z),M_n(\Z))$ on $M$ is a matrix in the $\Z$-algebra $\Z[M]$. Since $\Z[X]$ is obviously a subring of ${\rm Int}(M_n^p(\Z),M_n(\Z))$, we have:
\begin{equation*}
{\rm Int}(M_n^p(\Z),M_n(\Z))(M)=\{f(M)\,|\,f\in{\rm Int}(M_n^p(\Z),M_n(\Z))\}=\Z[M].
\end{equation*}
This result is a generalization of \cite[Theorem 6.4]{Frisch1} (indeed, it holds over any integral domain $D$). Let $C_p$ be the companion matrix of $p(X)$. By \cite[Lemma 2.2]{Per2} we have ${\rm Int}(M_n^p(\Z),M_n(\Z))={\rm Int}(\{C_p\},M_n(\Z))$. By the previous equality we have
\begin{equation*}
{\rm Int}(\{C_p\},M_n(\Z))={\rm Int}(\{C_p\},\Z[C_p])\doteqdot\{f\in \Q[X]\,|\,f(C_p)\in \Z[C_p]\}.
\end{equation*}
If we put together Corollary \ref{IntalphaZalpha=IntMnpalpha} and the previous equality we have
\begin{equation}\label{IntQalphaCpalpha}
{\rm Int}_{\Q}(\{\alpha\},\Z[\alpha])={\rm Int}(\{C_{p_\alpha}\},\Z[C_{p_\alpha}])
\end{equation}
In the same way, we can prove that ${\rm Int}_{\Q}(\Z[\alpha])={\rm Int}(\Z[C_{p_\alpha}])$.
\end{Rem}
\vskip0.3cm
\begin{Rem}\label{Mnirr}
By \cite[Proposition 6.2]{Frisch1} we have
\begin{equation}\label{SF}
\textnormal{Int}(M_n(\Z))=\textnormal{Int}(M_n^{\textnormal{irr}}(\Z),M_n(\Z)),
\end{equation}
where $M_n^{\textnormal{irr}}(\Z)$ is the set of matrices with irreducible characteristic polynomial. Let $\mathcal{P}_n^{{\rm irr}}$ be the subset of $\mathcal{P}_n$ of irreducible polynomials.
We can make a partition of $\mathcal{P}_n^{{\rm irr}}$ according to which number field of degree $n$  a polynomial $p\in\mathcal{P}_n^{\textnormal{irr}}$ has a root:
$$\mathcal P_n^{\textnormal{irr}}=\bigcup_{[K:\Q]=n} \mathcal P_n^K$$
where the union is taken over the set of all number fields $K$ of degree $n$ and $\mathcal P_n^K$ is the set of minimal polynomials of algebraic integers of $O_K$ of maximal degree $n$. Notice that for each $p\in \mathcal P_{n,K}$ we have $K\cong\Q[X]/(p(X))\cong\Q(\alpha)$, where $\alpha$ is a root of $p(X)$ (in particular, $\alpha$ is an algebraic integer of $K$). To ease the notation, we set $M_n^{\mathcal P_n^K}(\Z)=M_n^{K}(\Z)$. This is the set of matrices whose characteristic polynomial is irreducible and has a root in $K$. By (\ref{SF}) we have
\begin{equation}\label{IZMn}
{\rm Int}(M_n(\Z))=\bigcap_{[K:\Q]=n}{\rm Int}(M_n^{K}(\Z),M_n(\Z)).
\end{equation}
We have the following theorem, which resembles Theorem \ref{LW}:
\begin{Theorem}\label{3rdthm}
For a given number field $K$, the ring ${\rm Int}(M_n^K(\Z),M_n(\Z))$ is not integrally closed and its integral closure is ${\rm Int}_{\Q}(O_K)$.
\end{Theorem}
\vskip0.2cm
We prove Theorems \ref{2ndthm} and \ref{3rdthm} together in the next and final subsection. 
\end{Rem}
\vskip0.5cm
\subsection{\textbf{Proofs of Theorems \ref{2ndthm} and \ref{3rdthm}}}\label{ProofsThms}
\vskip0.4cm
The connection between rings of integral-valued polynomials and rings of integer-valued polynomials over matrices is given by Corollary \ref{IntalphaZalpha=IntMnpalpha} and (\ref{IntQalphaCpalpha}). Notice that the two rings $\Z[\alpha]$ and $\Z[C_{p_{\alpha}}]$ are isomorphic, since both are isomorphic to $\Z[X]/(p_{\alpha}(X))$.

We may now represent the rings of integer-valued polynomials over a set of matrices $M_n^P(\Z)$, $P$ a subset of $\mathcal{P}_n^{{\rm irr}}$, as an intersection of the pullback rings $R_\alpha$, where $\alpha$ ranges through the set of roots of the polynomials in $P$. Indeed, for such a subset $P$, by Corollary \ref{IntalphaZalpha=IntMnpalpha} and (\ref{critMnP}) we have:
$${\rm Int}(M_n^P(\Z),M_n(\Z))=\bigcap_{\alpha\in\mathcal{A}(P)}R_\alpha=\mathcal{R}_{\mathcal{A}(P)}$$
where $\mathcal{A}(P)\subseteq A_n$ is the set of roots of all the polynomials $p(X)$ in $P$ (notice that by assumption all these algebraic integers have degree equal to $n$). By Galois invariance (see Remark \ref{PropRASA}) for each polynomial $p\in P$ we can just take one of its roots. In this way we have
\begin{equation}\label{IntMnS}
{\rm Int}(M_n^P(\Z),M_n(\Z))=\mathcal{R}_{\mathcal{A}(P)}\subseteq \mathcal{S}_{\mathcal{A}(P)}={\rm Int}_{\Q}(\mathcal{A}(S),\mathcal{A}_n).
\end{equation}
By Theorem \ref{IntCloSARA}, the latter ring is the integral closure of ${\rm Int}(M_n^P(\Z),M_n(\Z))$. In particular, by Remark \ref{RAintclo}, ${\rm Int}(M_n^P(\Z),M_n(\Z))$ is integrally closed if and only if $\Z[\alpha]=O_{\Q(\alpha)}$,  for all $\alpha\in \mathcal{A}(S)$. 

Let $K$ be a number field of degree $n$ over $\Q$. Notice that the set of roots of the minimal polynomials of algebraic integers of degree $n$ in $K$ (this set is denoted by $\mathcal P_n^K$ in Remark \ref{Mnirr}) is exactly $O_{K,n}$, that is: $\mathcal{A}(\mathcal P_n^K)=O_{K,n}$. Then (\ref{IntMnS}) gives:
\begin{equation}\label{IntMnK}
{\rm Int}(M_n^K(\Z),M_n(\Z))=\mathcal{R}_{O_{K,n}}\subset\mathcal{S}_{O_{K,n}}={\rm Int}_{\Q}(O_K).
\end{equation}
The equality $\mathcal{S}_{O_{K,n}}={\rm Int}_{\Q}(O_K)$ is given by Theorem \ref{1stthm}. Hence, the integral closure of ${\rm Int}(M_n^K(\Z),M_n(\Z))$ is ${\rm Int}_{\Q}(O_K)$. This containment is also proper. In fact, the overrings $R_\alpha$, for $\alpha\in O_{K,n}$, are integrally closed if and only if $\Z[\alpha]=O_K$. In general there are plenty of $\alpha\in O_{K,n}$ such that this condition is not satisfied (see Proposition \ref{OnApoldenseO}). This concludes the proof of Theorem \ref{3rdthm}.

In the same way, by (\ref{IZMn}) and (\ref{IntMnK}) we have
\begin{equation}\label{IntMnZ}
{\rm Int}(M_n(\Z))=\bigcap_{[K:\Q]=n}{\rm Int}(M_n^{K}(\Z),M_n(\Z))=\bigcap_{[K:\Q]=n}\mathcal{R}_{O_{K,n}}=\mathcal{R}_{A_n}
\end{equation}
(the latter equality is given in Remark \ref{Annotpolclo}). We remark that, as already mentioned in \cite{LopWer}, this representation of ${\rm Int}(M_n(\Z))$ as an intersection of the rings $R_\alpha$ for $\alpha\in A_n$, shows that ${\rm Int}(M_n(\Z))$ is not Pr\"ufer, since there are many overrings $R_\alpha$ which are not integrally closed: by Proposition \ref{RSalphapullbacks} it is sufficient to consider an algebraic integer $\alpha$ of degree $n$ such that $\Z[\alpha]\subsetneq O_{\Q(\alpha)}$; then the corresponding $R_\alpha$ is such an overring. 

Finally, since for all $n$ we have ${\rm Int}(M_n(\Z))\subset{\rm Int}(M_{n-1}(\Z))$, by (\ref{IntMnZ})  $\mathcal{R}_{A_n}\subset\mathcal{R}_{A_{n-1}}$, so that 
$\mathcal{R}_{\mathcal A_n}=\bigcap_{m=1,\ldots,n}\mathcal{R}_{A_m}=\mathcal{R}_{A_n}$. In particular, by Theorem \ref{IntCloSARA}, the rings $\mathcal{S}_{A_n}$ and $\mathcal{S}_{\mathcal A_n}$ are equal, since they are the integral closure of $\mathcal{R}_{A_n}$ and $\mathcal{R}_{\mathcal A_n}$, respectively, and these two rings coincide. By (\ref{SAnAn}), the equality $\mathcal{S}_{A_n}=\mathcal{S}_{\mathcal A_n}$ is exactly Theorem \ref{2ndthm}.

\vskip0.8cm
\subsection*{\textbf{Acknowledgments}}
\noindent The author wishes to thank Stefania Gabelli for her survey article about pullbacks and Alan Loper for his useful suggestions.\\ 
The author was supported by the Austrian Science Foundation (FWF), Project Number P23245-N18.

\end{document}